\documentclass[12pt,thmsa]{amsart}
\usepackage{amsthm,amssymb,amstext}
\usepackage{amsfonts}
\usepackage{float}
\usepackage{amsmath}
\usepackage{latexsym}
\usepackage{amsopn}

\usepackage[cp1251]{inputenc}
\usepackage{graphicx}

\begin{document}
\thispagestyle{empty}
\title[On the geometry of the  domain of the solution of nonlinear Cauchy problem]
%On the geometry of the domain the solving nonlinear Cauchy problem]
{On the geometry of the domain of the solution of nonlinear Cauchy problem}
%On the geometry of the domain the solving nonlinear Cauchy problem}
\author[M. Z. Menteshashvili]{\'A. Figula (Debrecen) and M. Z. MENTESHASHVILI (Tbilisi)}
\address{\'A. Figula \\
Institute of Mathematics, University of Debrecen  \\
4010 Debrecen, P.O.Box 12. Hungary \\
and \\
M. Z. Menteshashvili\\
    Muskhelishvili Institute of Computational Mathematics\\
    of the Georgian Technical University\\
    Akuri str., 8, Tbilisi 0160\\
    Sokhumi State University, A. Politkovskaia str., 9\\
    Tbilisi 0186, Georgia.}

\email{figula@science.unideb.hu, and marimen1963@gmail.com}

\thanks { This paper are supported by the European Union's Seventh
Framework Programme (FP7/2007-2013) under grant agreements no. 317721, no. 318202.}

\begin{abstract}
The Cauchy and the inverse problems are considered  for a second order quasi-linear equation with an admissible parabolic degeneration on open and closed data support.
The inverse problem is investigated in the cases that the definition domain of the solution of the differential equation contains a  gap.

\end{abstract}
\maketitle

\bigskip
\noindent 
\centerline{\bf 1. Introduction}

\bigskip 
In this paper we consider the following Cauchy problem for second order hyperbolic differential equations: find a solution 
$u(x,y)$ of the equation by the initial conditions  $u|_{y=0}=\tau(x)$, $u_y|_{y=0}=\nu(x)$, where
$\nu(x) \in C^1(\mathbb R)$,
$\tau(x) \in C^2(\mathbb R)$ are given functions such that  $\nu(x)$ is once-, and  $\tau(x)$ is twice-continuously differentiable.
To solve this problem one can use the method of characteristics (cf. \cite{9}, pp. 164-166). The families of characteristic curves for linear differential equations are defined by the principal part of the equation 
(cf. \cite{B2}, p. 3). The class of hyperbolic equations is defined through the characteristic roots by the inequality 
$\lambda_1\neq\lambda_2$. Characteristic directions are defined at every point by the relations
 $\frac{dy}{dx}=\lambda_1(x,y,u,u_x,u_y)$, $\frac{dy}{dx}=\lambda_2(x,y,u,u_x,u_y)$ (cf. \cite{M}).
	
As soon as the Cauchy problem is stated for a hyperbolic equation, it is needed to impose certain conditions on the initial data support. It is required that every characteristic curve of the equation intersects the initial support at most one point. Furthermore, it is required  that there are no points of the support at which the support has the characteristic direction (cf. \cite{9}, p. 166). If equations are considered along the characteristic manifolds, then we see that these properties are generated by the equations themselves.

For clarity, let us consider the string vibration equation $u_{xx}-u_{yy}=0$ on the plane of variables
$x$, $y$. The one-parameter families of characteristic curves are related to the characteristic roots $\lambda_1=1$ and 
$\lambda_2=-1$ and expressed by the relations $x-y={\it \hbox{const}}$, $x+y={\it \hbox{const}}$. The combination of the first derivatives of the solution $u_x-u_y$ remains constant along any characteristic curve of the first family, whereas the combination $u_x+u_y$ remains constant along any characteristic curve of the second family. These combinations are called the characteristic invariants (cf. \cite{12}). If the initial support $y=f(x)$ of the Cauchy problem intersects a characteristic curve of the first family at two points $(x_0,y_0)$ and $(x_1,y_1)$, then using the above property of the characteristic invariants we obtain that $u_x(x_0,y_0)-u_y(x_0,y_0)=u_x(x_1,y_1)-u_y(x_1,y_1)$. This implies that the derivatives of the solution are bounded and therefore stating the problem the initial data on the support cannot be chosen arbitrary and thus they are redundant to make the problem well-posed. Hence the conditions for the initial data and the initial support is quite natural. According to Hadamard (cf. \cite{H}) the Cauchy problem can appear to be well-posed in some cases involving the closed initial support. Such problems were studied by Aleksandryan \cite{A}, Sobolev \cite{S},  Vakhania \cite{V}, Wolfersdorf \cite{W}  and others for linear equations.

The above facts are true in the case of linear equations. However we have a different situation if the families of characteristic curves cannot be defined a priori since they depend on the unknown solution $u$ and its derivatives $u_x$, $u_y$. In particular this  situation  takes place for the non-linear hyperbolic equation
\begin{equation} \label{equ1} a u_{xx}+b u_{xy}+c u_{yy}=f, \end{equation}
where the principal coefficients $a, b, c$ and the right-hand side of (\ref{equ1}) are functions of five variables $x$, $y$, $u$, $u_x$, $u_y$ (cf. \cite{B1}).  Such cases are somewhat difficult to investigate since they demand modification and generalization of linear problems (cf. \cite{Gv}). Among the papers based on the application of the method of characteristics to non-linear hyperbolic problems we refer to \cite{7}, \cite{GMB}, \cite{4}, \cite{Gv1},  \cite{Men}, where the structure of the definition domains of the solutions and that of influence domains of initial and characteristic perturbations are investigated in singular cases. Another type of the application of the method of characteristics to non-linear equations is the initial-characteristic Darboux problem. The formulations of Darboux problems for non-linear  equations should be revised taking into account general characteristic invariants (cf. \cite{7}, \cite{BM}, \cite{3}, 
\cite{4}, \cite{6}, \cite{5}). In \cite{BM} the authors formulated  correctly an initial-characteristic Darboux problem for the quasi-linear equation
\begin{equation}\label{e1}
    x^2(u_y^4u_{xx}-u_{yy})=cuu_y^4, \;\; c=const,
\end{equation}
which arises by study nonlinear oscillations.

The aim of our investigation is to discuss some questions stating an initial problem on the data support
$[a,b] \in \mathbb R$ for the following second order non-strictly hyperbolic equation
\begin{equation} \label{equ2} (u^2_y-1)u_{xx}-2u_xu_yu_{xy}+u^2_xu_{yy}=0   \end{equation}
which is called Dubreil-Jacotin equation (cf. \cite{B1}, p. 442).   We consider also the inverse problem for equation 
(\ref{equ2}) on some open and closed data supports. Equation (\ref{equ2}) is interesting for physical applications since the 
two-dimensional flow of an inviscid incompressible fluid in gravitational field can be decribed by it (cf. \cite{Yu}, 
p. 535).

The characteristic roots of (\ref{equ2}) are
$$-\frac{u_x}{u_y+1},\; -\frac{u_x}{u_y-1}.$$	
Therefore equation (\ref{equ2}) is hyperbolic everywhere, except for the points at which the derivative of the sought solution  $u_x$ has zero values. At these points, equation (\ref{equ2}) parabolically degenerates. Hence equation 
(\ref{equ2}) has mixed type
(cf. \cite{B2}). It is one of the most important problems of mathematical physics to study the properties of solutions of equations of mixed type (cf. \cite{2}, \cite{T}, \cite{Se}, \cite{10}). The first fundamental results in this direction were obtained by the Italian mathematicians Francesco Tricomi in the twenties of the $19$th century.

The Cauchy problem of equation (\ref{equ2}) such that the given functions $\tau $, $\nu $ are defined in the circumference of the unit circle, i.e. initial problem on the closed data support $\gamma:\,\rho=1,$ $x=\rho \cos\vartheta,\,y=\rho \sin\vartheta$, $0 \le \vartheta \le 2 \pi$, was investigated in \cite{Men}. 

In this paper we study the Cauchy problem of equation (\ref{equ2}) on the support $[a,b]$ (see Problem 1 in section 2). Since the set of points of parabolic degeneration of the considered equation is not defined a priori, formulating the Cauchy problem  it is always necessary to ascertain whether the equation degenerates on the data support or not. This rule should also be observed in the case to be considered below. 
Suppose that the initial support is given in the explicit form by the equation $y=f(x)$, where the function $f$ satisfies certain conditions which will be discussed below. Then we can define the direction of any characteristic curve that comes out from an arbitrary point $x_0$ of the support by the relation 
\begin{equation} \frac{dy}{dx}= \lambda_1 \left(\tau(x_0), \frac{\tau'(x_0)- \nu(x_0) f'(x_0)}{1+(f'(x_0))^2}, 
\frac{\nu(x_0)+ \tau'(x_0) f'(x_0)}{1+(f'(x_0))^2} \right). \nonumber \end{equation} 
Taking $x_0$ as a parameter we can define the direction of the characteristic curves at every point of the support. Though the situation for the second family of characteristic curves is treated analogously, to solve the Cauchy and inverse problems for equation (\ref{equ2}) we have to write equations for the characteristic curves of both families. To establish this equations we use the representation formula of the general integral for (\ref{equ2}). After this we consider the set 
of intersection points for both families of characteristic curves. This set creates the domain within which the initial support is located. In Theorem 1 we give the solution of the Cauchy problem of (\ref{equ2}) in integral form. 

The combinations of first derivatives of the solution which are constant for the string vibration equation differ from those for non-linear equations (\ref{equ1}). In the latter case these combinations depend on an unknown solution and its first derivatives. In aerodynamics, these combinations are called  Riemann invariants (cf. \cite{B2}, p. 23).
It is rather difficult to use the characteristic invariants for the solution of problems in the case of quasi-linear equations. The characteristic invariants for equation $(\ref{equ2})$ are $u+y= const$ and $u-y=const$. Using these invariants we prove that the considered inverse problem (see Problem 2 in section 3) has a solution under certain  regularity condition (cf. Theorem 2).  We illustrate the Cauchy and the inverse problems in Examples 1, 2, 3. In our examples we deal with the interesting cases that the families of the characteristic curves have either common envelopes or singular points.

\bigskip
\noindent
\centerline{\bf 2. Cauchy problem on equation (\ref{equ2})}

\bigskip 
{\bf Problem 1.} {\it Let $\tau(x)$, respectively $\nu(x)$
be real functions, which are
 twice, respectively once continuously differentiable. Find a function  $u(x,y)$,
which satisfies equation (\ref{equ2}) and the conditions
\begin{equation} \label{equ3} u(x,0)=\tau(x), \quad        u_y(x,0)=\nu(x), \quad        a \le x \le b. \end{equation}
It is also required to find a domain where the solution can be completely defined.}
\vskip 0.3cm
{\sc\bf Theorem 1.} {\it If $\tau'(x) \neq 0$ for all $x \in [a,b]$, then the integral of the problem (\ref{equ2}),
(\ref{equ3}) can be  written into the form
\begin{equation} \label{equ0} x= a +\frac{1}{2} \int^{T(u+y)}_a (1-\nu(t)) dt + \frac{1}{2} \int^{T(u-y)}_a   (1+\nu(t)) dt, \end{equation}
where
%$$F(t)= , \quad G(t)= \frac{1+\nu(t)}{2\tau^\prime(t)},$$
%and
$x=T(z)$ denotes the inverse function of $z=\tau(x)$. }
\vskip 0.3cm
\begin{proof}
We start the investigation of the problem by considering the form of the general integral of equation (\ref{equ2}). It can be represented by the following functional equation
\begin{equation} \label{equ5} f(u+y)+g(u-y)=x, \nonumber \end{equation}
where $f$, $g$ are arbitrary functions (see \cite{Gv}).
In order to provide the required smoothness of the solution $u(x,y)$, here it is assumed that the arbitrary functions $f$,
$g$ belong to the class $C^2(\mathbb R)$. Considering the formula
$f(u(x,0)+0)+g(u(x,0)-0)=x$ and taking the derivation of it with respect to the variables $x$ and $y$
we obtain
\begin{equation} \label{equ6} f'(u(x,0)+0) u_x(x,0)+ g'(u(x,0)-0) u_x(x,0)=1 \end{equation}
\begin{equation} \label{equ7} f'(u(x,0)+0)(u_y(x,0)+1)+g'(u(x,0)-0)(u_y(x,0)-1)=0. \end{equation}
Hence we come to a system of two linear algebraic equations with respect to the variables $f^\prime(u(x,0))$, $g^\prime(u(x,0))$. Solving this system of linear equations one gets
\begin{equation} \label{equhatuj} f'(u(x,0))=\frac{1-\nu(x)}{2\tau^\prime(x)} :=F(x) \end{equation}
\begin{equation} \label{equhetuj} g'(u(x,0))=\frac{1+\nu(x)}{2\tau^\prime(x)} :=G(x). \end{equation}
If the condition $\tau^\prime(x) \neq 0$ is fulfilled for all $x \in [a,b]$, then on the closed interval
$[a,b]$ the equation $\tau(x)= z$ is uniquely solvable in the class of real solutions. We denote this solution by
$x=T(z)$.
Integrating relations (\ref{equhatuj}), (\ref{equhetuj}) between $a$ and an arbitrary value $x=T(z)$ we obtain
\begin{equation} \label{equ8}
f(z)-f(\tau(a))= \int \limits_{\tau(a)}^z f'(v) dv=\int \limits_{a}^x F(t) \tau'(t) dt, \nonumber \end{equation}
\begin{equation} \label{equ8uj}
g(z)-g(\tau(a))= \int \limits_{\tau(a)}^z g'(v) dv=\int \limits_{a}^x G(t) \tau'(t) dt, \nonumber \end{equation}
or equivalently
\begin{equation} \label{equ9}
f(z)=f(\tau(a))+\int \limits_{a}^{T(z)} F(t) \tau'(t) dt, \end{equation}
\begin{equation} \label{equ10}
g(z)=g(\tau(a))+\int \limits_{a}^{T(z)} G(t) \tau'(t) dt. \end{equation}
Their sum already yields an implicit solution of problem (\ref{equ2}), (\ref{equ3}) which contains the undefined free functions for
$x=a$. They can be identified by normalization. In particular, if for the functions $f$, $g$ defined by
(\ref{equ9}), (\ref{equ10}), relation (\ref{equ5}) is fulfilled at the point $(a,0)$, then we have
\begin{equation} \label{equ11}
f(\tau(a))+ g(\tau(a)) + \int^{T[u(a,0)+0]}_a   F(t) \tau'(t) dt + \int^{T[u(a,0)-0]}_a   G(t) \tau'(t) dt = a.
\end{equation}
	Note that the identities $T(\tau(x))=x$
are valid for all values of $x \in[a,b]$, including $x=a$, and therefore the upper bounds of both integrals in (\ref{equ11}) are $a$. Hence one has $f(\tau(a)) + g(\tau(a))= a.$
Finally, the implicit solution of problem (\ref{equ2}), (\ref{equ3}) can be written as
\begin{equation} \label{equ12} x-a =\int^{T(u+y)}_a  F(t) \tau'(t) dt + \int^{T(u-y)}_a  G(t) \tau'(t) dt, \end{equation}
where the functions $F(t)$ and $G(t)$ are defined by (\ref{equhatuj}), (\ref{equhetuj}).
This representation is obtained if condition $\tau'(x) \neq 0$ is fulfilled for all $x \in [a,b]$ and the expressions under the integral sign are integrable.
\end{proof}
\vskip 0.3cm
In order to establish the definition domain of solution (\ref{equ12}), it is necessary to investigate the structure of the characteristic curves. Keeping in mind representation (\ref{equ12}) and that the relation $u+y=const$ is fulfilled along the characteristic curves of the first family, we can obtain an equation for each characteristic of this family that passes through an arbitrary point $x$ of the data support $[a,b]$. At this arbitrarily chosen point, with which we associate the argument
$x_0,$ we have
$$u(x_0,\,0) = \tau(x_0). $$
Hence the relation
\begin{equation} \label{rel1} u(x,y)+ y = \tau(x_0) \end{equation}
 is fulfilled along the characteristic curve satisfying the relation $u+y=const$ and passing through the point  
$(x_0,\,0)$. Since the right-hand side of this equality is a completely defined expression, we conclude that
\begin{equation} \label{rela2} u(x,y)- y = \tau(x_0) -2y. \end{equation} 
Substituting (\ref{rel1}) and (\ref{rela2}) into (\ref{equ12}), we obtain the equation for a characteristic curve of the first family in the explicit form
\begin{equation} \label{exp1} x-a= \int^{T(\tau(x_0))}_a  F(t) \tau^\prime (t) dt + \int^{T(\tau(x_0)-2y)}_a G(t) \tau^\prime(t) dt. \end{equation}
Keeping in mind that the identity $T(\tau(x_0))=x_0$ is valid, the upper bound of the first integral in (\ref{exp1}) will be $x_0$. If we introduce the notation $x_0= c$, then this expression becomes the parameter which takes values from the interval $[a,b]$.  Therefore, all characteristic curves of the family $u+y=const$ which pass through the points of the data support have the representation
\begin{equation} \label{equ13} x= a + \int^{c}_a  F_1 (t) dt + \int^{T(\tau(c)-2y)}_a   G_1 (t) dt, \end{equation}
where
\begin{equation} \label{firstchar} F_1(t)= F(t) \tau^\prime(t)=\frac{1- \nu(t)}{2} \nonumber \end{equation}
\begin{equation} \label{firstchar2} G_1(t)= G(t) \tau^\prime(t)=\frac{1+\nu(t)}{2}. \nonumber \end{equation}

Analogously, for all characteristic curves of the family $u-y=const$ we obtain

\begin{equation} \label{equ14} x=a+ \int^{T(\tau(c)+2y)}_a   F_1(t)dt +  \int^{c}_a   G_1(t)dt. \end{equation}
Here the parameter $c$ takes its values from the interval $[a,b]$.  For both families, in (\ref{equ13}) and (\ref{equ14}) there is a one-to-one correspondence between the equation for a characteristic curve and the parameter value. This fact is stipulated by the condition
$\tau'(x) \neq 0$ for all $x \in [a,b]$.

\bigskip
\noindent
\centerline{\bf 3.\,\,Inverse problems}

\bigskip
\noindent Now we consider certain inverse problems of the Cauchy problem for the Dubreil-Jacotin equation (\ref{equ2}).
Below we will treat a variant of inverse problem which require the construction of a given equation under the condition that the characteristic curves of both families are known a priori. The cases that the families of characteristic curves have either common envelopes or singular points  are particularly interesting to investigate. In these situations, gaps may appear in the definition domain of a solution  of the problem.  As has been mentioned above, equation (\ref{equ2}) is of hyperbolic type and the set of parabolic degeneration is not defined a priori because it depends on the behavior of an arbitrary solution $u$ and with respect to the variable $x$. The characteristic invariants (cf. \cite{B2}, p. 23) corresponding to equation (\ref{equ2}) are given as follows

\begin{equation} \label{firstcase}  u+y=const, \quad {\text{\it for}} \quad \lambda_1=\frac{u_x}{u_y+1},  \end{equation}
\begin{equation} \label{secondcase} u-y=const, \quad {\text{\it for}} \quad \lambda_2=\frac{u_x}{u_y-1}. \end{equation}
From these equations we conclude that the families of characteristic curves of the equation are not defined a priori because they depend on the sought solution and on the first order derivatives. In the case of equation (\ref{equ2}), we can admit in the role of characteristics the families of characteristic curves along which relations (\ref{firstcase}) and (\ref{secondcase}) are fulfilled. Characteristic families can be therefore given a priori in an arbitrary manner. The statement of the inverse problems considered here rests on this fact.

Let us consider two one-parameter families of plane curves which are given explicitly by the equations
\begin{equation} \label{curve1} y=\varphi_1(x,c), \end{equation}
\begin{equation} \label{curve2} y=\varphi_2(x,c), \end{equation}
where $\varphi_1$, $\varphi_2$  are given, twice differentiable functions with respect to the variable $x$. Let $\varphi_1$, $\varphi_2$ be defined for any value of the real parameter $c$. Assume that any curve of these families necessarily intersects once the straight line $y=0$. We denote by $D_1$ the domain of the plane $(x,y)$ which is completely covered by the family of characteristic curves given by equation (\ref{curve1}) if the parameter $c$ runs continuously through all real values. Analogously, we denote by $D_2$ the domain  of the plane $(x,y)$ which is completely covered by the characteristic lines given by equation (\ref{curve2}). Also, we introduce the notation $D=D_1\cap D_2$,  $I=D\cap \{y=0\}.$
\vskip 0.3cm
{\sc\bf Problem 2 (Inverse problem).}
{\it Find the initial conditions of a regular solution $u(x,y)$ of equation (\ref{equ2}) and its derivative with respect to the normal direction on the interval $I$ of the straight line $y=0$ if the families of plane curves given by equations (\ref{curve1}) are the characteristics corresponding to invariants (\ref{firstcase}), while (\ref{curve2}) is the family of characteristics corresponding to invariants (\ref{secondcase}).}
\vskip 0.3cm
{\sc\bf Theorem 2.} {\it If the condition $\varphi_1^\prime (x,c^\prime) \neq \varphi_2^\prime (x,c^{\prime\prime})$, is fulfilled for all $c^\prime$, $c^{\prime\prime}\in \mathbb R$, then there exists a solution of Problem 2.}
\vskip 0.3cm
\begin{proof}  To investigate the posed problem, we need a structural analysis of the characteristic curves of equation (2). The characteristic roots corresponding to equation (2) define at every point two characteristic directions
\begin{equation} \label{directionfirst} \frac{dy}{dx}= -\frac{u_x}{u_y+1}\equiv\lambda_1(x,y),\end{equation}
\begin{equation} \label{directionsecond} \frac{dy}{dx}= -\frac{u_x}{u_y-1}\equiv\lambda_2(x,y).\end{equation}
It follows from the posed problem that the family of characteristic curves defined by equation (\ref{curve1}) corresponds to the root $\lambda_1$ given by (\ref{directionfirst}), while that given by equation (\ref{curve2}) corresponds to the root $\lambda_2$ given by (\ref{directionsecond}).  Therefore we can calculate the values of first derivatives $u_x$, $u_y$ of the unknown solution $u(x,y)$ along any characteristic curve. Note that the parameters contained in families (\ref{curve1}) and (\ref{curve2}) can be defined by the abscissa  of the intersection point of the concrete curve and the axis $y=0$. Indeed, let the curve corresponding to the parameter  $c^*,$ intersect the straight line $y=0$ at a point $x_0$
$$\varphi_1(x_0,c^*)=0.$$
Solving the equation for the parameter $c^*(x_0)=c,$ family (\ref{curve1}) takes the form
$$y=\varphi_1(x,x_0),$$
which is equal to zero for $x=x_0$. If family (\ref{curve1}) corresponds to the characteristic root $\lambda_1$, then the relation
\begin{equation} \label{lambda1} \frac{d\varphi_1(x,x_0)}{dx}= -\frac{u_x(x,\varphi_1(x,x_0))}{u_y(x,\varphi_1(x,x_0))+1}.\end{equation}
is fulfilled. Analogously, along the second family of the characteristic curves the following equality is fulfilled
\begin{equation} \label{lambda2} \frac{d\varphi_2(x,x_0)}{dx}= -\frac{u_x(x,\varphi_2(x,x_0))}{u_y(x,\varphi_2(x,x_0))-1}.\end{equation}
Equalities (\ref{lambda1}), (\ref{lambda2}) have the following at the points of the interval $I$:
\begin{equation} \label{lambda1uj}\frac{d\varphi_1(x,x_0)}{dx}\mid_{x=x_0}= -\frac{\tau^\prime(x_0)}{\nu(x_0)+1}, \;\,x_0\in I,\end{equation}
\begin{equation} \label{lambda2uj}\frac{d\varphi_2(x,x_0)}{dx}\mid_{x=x_0}= -\frac{\tau^\prime(x_0)}{\nu(x_0)-1}, \;\,x_0\in I.\end{equation}
Since the point $(x_0,0)$ is arbitrarily chosen from the interval $I$, $x_0$ in equations (\ref{lambda1uj}), (\ref{lambda2uj}) can be replaced with $x$. Hence we can easily define the unknown functions $\tau^\prime$ and $\nu$. Integrating the function $\tau^\prime$, we finally obtain the solution $\tau(x)$, $\nu(x)$ of the problem. Naturally, the function $\tau$ is defined up to an integration constant which will be defined uniquely if we give the value of the function $u(x,y)$ at one arbitrary point of the interval $I$. The theorem is thereby proved.
\end{proof} 

\bigskip
\noindent 
From the above general argumentation we can draw a conclusion. If families (\ref{curve1}) and (\ref{curve2}) do not have common characteristic directions at anyone of the points, i.e. 
if equation (\ref{equ2}) is of strictly hyperbolic type, then the solution of the inverse problem presents no difficulty 
(see Example 1). However, if the equation parabolically degenerates on certain set of points, then the situation drastically changes. 
We illustrate this in Examples 2, 3.  

\bigskip
\noindent 
{\sc\bf Example 1.}
We solve the Cauchy- and inverse problems in the following concrete example.
\newline
\noindent 
Cauchy problem:
Let $u(x,0)=\tau(x)$ be the function $x$ and $u_y(x,0)=\nu(x)$ be the function
$1-e^x$. According to (\ref{equ0}) the solution (\ref{equ12}) of the Cauchy problem in implicit form is given by equation
$$\frac{1}{2}e^{u+y}+u-y- \frac{1}{2} e^{u-y}-x=0.$$ 
The characteristic curves (\ref{equ13}) of the family $u+y=const$, as parametric curve $f(x,y,c)$ with respect to the parameter $c$ are given by
\begin{equation} \label{equuj15} 0=\frac{1}{2} e^c+c-2y- \frac{1}{2}e^{c-2y} -x=f_1(x,y,c) \end{equation}
The characteristic curves (\ref{equ14}) of the family $u-y=const$, as parametric curve $f(x,y,c)$ with respect to the parameter $c$ are given by
\begin{equation} \label{equuj16} 0=\frac{1}{2} e^{c+2y}+c- \frac{1}{2}e^{c} -x=f_2(x,y,c). \end{equation}
Differentiating the function $f_1(x,y,c)$ in (\ref{equuj15}) with respect to the parameter $c$ we get
\begin{equation} \label{equuj17} 0=\frac{1}{2} e^c (1-e^{-2y})+1=f_{1,c}(x,y,c). \end{equation}
Expressing from (\ref{equuj17}) the parameter $c$ we have
\begin{equation} \label{equuj18} c=\ln \left( \frac{2 e^{2y}}{1-e^{2y}} \right).  \end{equation}
Putting expression (\ref{equuj18}) into equation (\ref{equuj15}) we obtain 
\begin{equation} \label{equuj19} e^{1+x}= \frac{2}{1-e^{2y}}.  \end{equation}
Since the derivative $\frac{\partial f_{1,c}(x,y,c)}{\partial y}= e^{c-2y}$ is non-zero and the derivative
$\frac{\partial f_{1,c}(x,y,c)}{\partial c}= \frac{1}{2}(1- e^{-2y}) e^c$ is non-zero for all $y \neq 0$ we obtain that the envelope of the characteristic curves (\ref{equuj15}) of the first family is given by (\ref{equuj19}). An analogously computation shows that this envelope is also the envelope of characteristic curves (\ref{equuj16}) of the second family.
In the plane $\{ x, y \}$ the domain which is below of envelope (\ref{equuj19}) there does not exist solution of the Cauchy problem (see Figure 1 for parameter $c=-1$).

\begin{figure}[hptb]
\centering
 \includegraphics[width=11cm]{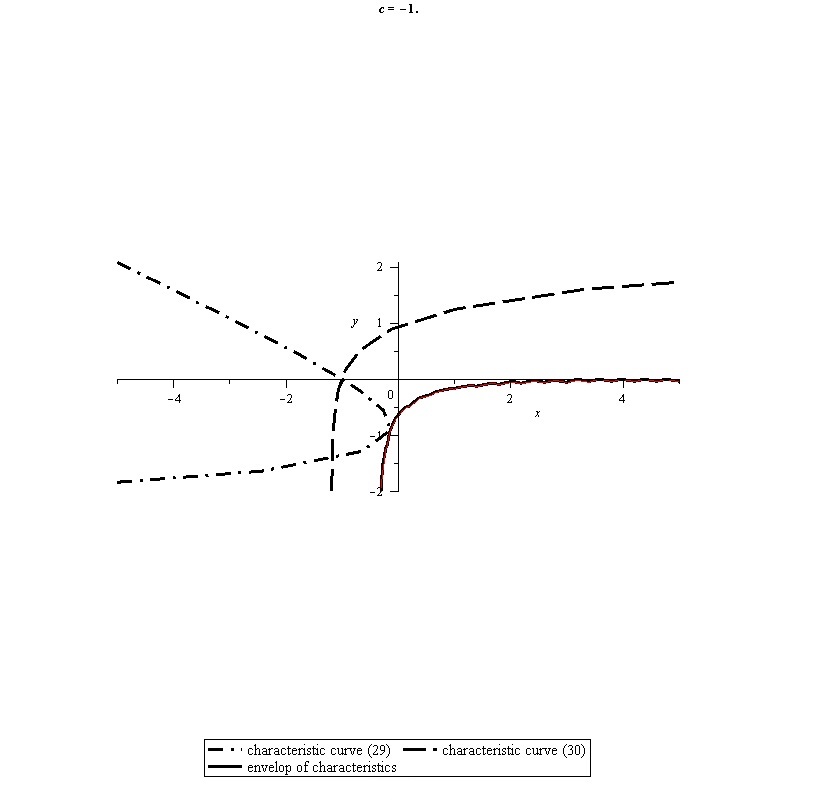}

 \caption{The characteristic curves and their envelope.}
\end{figure}

\noindent
Inverse problem:
Assume that the families of characteristic curves are given by:
\begin{equation} \label{curf1} x=\varphi_1(y,c) \equiv {\frac{1}{2}e^c-\frac{1}{2}e^{c-2y}+c-2y}, \quad {\text{\it for}} \quad \lambda_1, \end{equation}
\begin{equation} \label{curf2} x=\varphi_2(y,c) \equiv {-\frac{1}{2}e^c+\frac{1}{2}e^{c+2y}+c}, \quad {\text{\it for}} \quad \lambda_2. \end{equation}
From (\ref{lambda1uj}), (\ref{lambda2uj}) we get
\begin{equation} \label{exe1} -\frac{u_x(x_0,0)}{u_y(x_0,0)+1} = \frac{1}{e^{x_0}-2}, \end{equation}
for family (\ref{curf1}) and 
\begin{equation} \label{exe2} -\frac{u_x(x_0,0)}{u_y(x_0,0)-1} = \frac{1}{e^{x_0}}, \end{equation}
for family (\ref{curf2}).

The right-hand sides of the equalities (\ref{exe1}), (\ref{exe2}) depend only on $x_0$. From the relations (\ref{exe1}) and (\ref{exe2}) for the point $(x_0,0)$ we find the values of the functions
$\tau^\prime$  and $\nu$:
\begin{equation} \label{exe7} \tau^\prime (x) =  1,\end{equation}
\begin{equation} \label{exe8} \nu(x)=1-e^{x}.\end{equation}
From equality (\ref{exe7}) we obtain by integration the value of the solution $u(x,y)$ on the axis $y=0$:
\begin{equation} \label{exe9} \tau(x) = x+c,\;\;c=const.\end{equation}
The constant $c$ in (\ref{exe9}) can be obtained if we know the value of $\tau(x)$ in one arbitrary point.

\bigskip
\noindent 
{\sc\bf Example 2.} Now we consider the following inverse problem. 
Assume  the families of characteristic curves have common envelopes and are given in the form
\begin{equation} \label{curve3} x=\varphi_1(y,c) \equiv c - a\sqrt{\frac{a}{y+b}- 1},\,\;  -b\leq y \leq a-b, \end{equation}
\begin{equation} \label{curve4} x=\varphi_2(y,c) \equiv c + a\sqrt{\frac{a}{y+b}- 1},\,\; -b \leq y \leq a-b, \end{equation}
\begin{equation} a > b > 0, \; c= const. \nonumber \end{equation}
Each curve of family (\ref{curve3}) defined on an interval $x\in(-\infty,c]$ monotonically increases, has the asymptote 
$y=-b$. For
$x=c$ each curve of family (\ref{curve3}) is tangent to the straight line $y=a-b$ and therefore is tangent to one of the curves of family (\ref{curve3}). One of the curves of family (\ref{curve3}) passes at every point of the straight line
$y=a-b$ and smoothly continues from the same point to the completely defined curve of family (\ref{curve4}). Therefore both characteristic directions at the points of this straight line coincide. Hence the straight line $y=a-b$ is the line of parabolic degeneration for equation (\ref{equ2}). The straight line $y=-b$ is also the line of parabolic degeneration because all characteristic curves of both families are tangent to the straight line at infinity.

By direct calculations we establish that 
\begin{equation} \label{equuj1} \frac{u_x(x,y)}{u_y(x,y)-1} = \frac{2a^3(x-c_1)}{\left[(x-c_1)^2 + a^2 \right]^2}, \,\,\,\,  c_1 = x_0 + a\sqrt{\frac{a-b}{b}}, \end{equation}
holds for family (\ref{curve3}) and that  
\begin{equation} \label{equuj2} \frac{u_x(x,y)}{u_y(x,y)+1} = \frac{2a^3 (x-c_2)}{\left[(x-c_2)^2 +a^2\right]^2},\,\,\,\,   c_2  = x_0 - a\sqrt{\frac{a-b}{b}}, \end{equation}
is satisfied for family (\ref{curve4}). Note that the right-hand sides of the equalities (\ref{equuj1}), (\ref{equuj2}) are represented by expressions depending solely on the variable $x$. Relations (\ref{equuj1}) and (\ref{equuj2}) for the point 
$(x_0,0)$ yield the equalities
\begin{equation} \label{equuj3} \frac{\tau^\prime (x_0)}{\nu(x_0)-1} = \frac{2a^3 (x_0 -c_1)}{[ (x_0 -c_1 )^2 +a^2 ]^2}, \nonumber \end{equation}
for the first family of the characteristic curves, and
\begin{equation} \label{equuj4} \frac{\tau^\prime (x_0 )}{\nu(x_0)+1} = \frac{2a^3 (x_0 -c_2)}{[ (x_0 -c_2)^2  +a^2 ]^2},
\nonumber \end{equation}
for the second family.
Let us assume that the point $x_0\in(-\infty,+\infty)$ is arbitrary, then for any $(x,0)$ we obtain
\begin{equation} \label{equuj5} \frac{\tau^\prime (x)}{\nu(x)+1} = -\frac{2b^2}{a^2}\sqrt{\frac{a-b}{b}},\end{equation}
\begin{equation} \label{equuj6} \frac{\tau^\prime(x)}{\nu(x)-1} = \frac{2b^2}{a^2}\sqrt{\frac{a-b}{b}}.\end{equation}
Hence we find the values of the functions $\tau^\prime$  and $\nu$:
\begin{equation} \label{equuj7} \tau^\prime (x) =  -\frac{2b^2}{a^2} \sqrt{\frac{a-b}{b}},\end{equation}
\begin{equation} \label{equuj8} \nu(x)=0.\end{equation}
From equality (\ref{equuj7}) we find by integration the values of the solution $u(x,y)$ on the axis $y=0$:
\begin{equation} \label{equuj9} \tau(x) = -\frac{2b^2}{a^2}\sqrt{\frac{a-b}{b}}\,x+c,\;\;c=const.\end{equation}

As relations (\ref{equuj5}), (\ref{equuj6}) show the characteristic invariants $\frac{u_x}{u_y-1},\,\frac{u_x}{u_y+1}$ of both families take  a constant value. Despite this, we cannot assert that the straight line $y=0$ is the line of parabolic degeneration. To make such an assertion it is necessary that the invariants $u+y$ and $u-y$ are constant along $y=0$. The constancy of the invariants is due to the fact that families (\ref{curve3}), (\ref{curve4}) are sets of curves obtained by a parallel transfer along the axis of abscissas. Therefore every curve of both families has one and the same slope with respect to the axis of abscissas. This indicates that the values of the derivative $\frac{dy}{dx}$  along this axis preserve constancy.

\begin{figure}[hptb]
\centering
 \includegraphics[width=11cm]{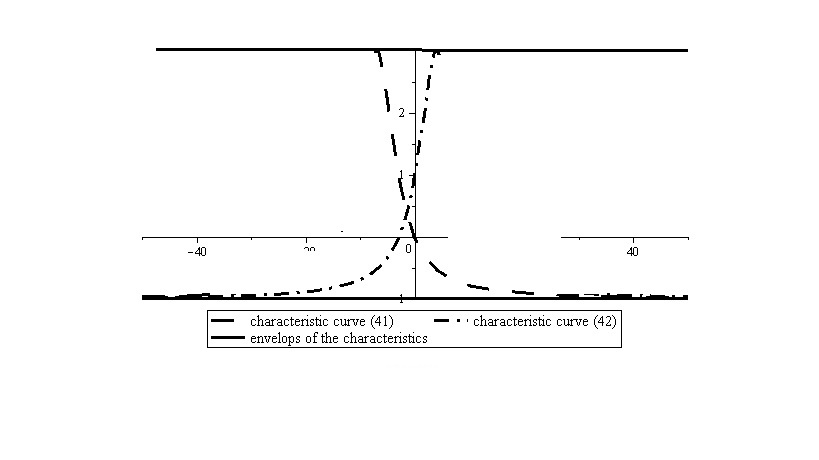}

 \caption{The characteristic curves and their envelopes.}
\end{figure}

\bigskip
\noindent
{\sc\bf Example 3.} 
The situation is more difficult if the families of the characteristic curves have, besides the lines of parabolic degeneration, also common singular points. As an example let us consider the case that families (\ref{curve1}),
(\ref{curve2}) have the common node. Each contour of the one-parameter family of curves
\begin{equation} \label{kor1} F(r,\vartheta,\vartheta_0) = \frac{2b(1-r\cos\vartheta)}{r^2-2r\cos\vartheta+1} 
-a^3\left\{ a^2 + \left[\frac{2br\sin\vartheta}{r^2-2r\cos\vartheta +1}+\vartheta_0\right]^2\right\}^{-1} = 0,\;\;
\end{equation}
$0\leq\vartheta \leq 2\pi$ is closed.  The whole family lies in the half plane to the left of the straight line $x =1$. The circumference
\begin{equation} \label{kor2} [(x-1)a + b]^2  + y^2a^2 = b^2 \end{equation}
is a common envelope of a family of curves.
The point $(1,0)$ is the node of the considered family. If $a>b$, then the unit circumference completely lies in the definition domain of family (\ref{kor1}). Every curve of this family  intersects twice the unit circumference. Every contour of family (\ref{kor1}) can be represented as the union of two arcs. The first arc is considered from the point
$(1,0)$ to the point of tangency to circumference (\ref{kor2}) if the movement occurs in the positive direction. The remaining part is considered to be the arc of the second family. If the polar angle of the point of tangency of the concrete curve to circumference (\ref{kor2}) is denoted by $\vartheta^*$, then family (\ref{kor1}) can be divided into two parts:
\begin{equation} \label{kor3} F_1(r,\vartheta,\vartheta_0) = 0,\;\;0\leq\vartheta \leq \vartheta^*, \end{equation}
\begin{equation} \label{kor4} F_2(r,\vartheta,\vartheta_0) = 0,\;\;\vartheta^* \leq \vartheta \leq 2\pi. \end{equation}
\vskip 0.3cm
{\sc\bf Problem 3 (Inverse problem).}
{\it Find the values of a regular solution  of equation (\ref{equ2}) and its derivative with respect to the normal direction on the circumference $x^2+y^2=1$, i.e. find the functions 
\begin{equation} \label{kor5} \lim\limits_{r\to 1} u(r,\vartheta)= \tau(\vartheta), 
\qquad \lim\limits_{r\to 1} u_r(r,\vartheta)= \nu(\vartheta),
\end{equation}
if the family of plane curves (\ref{kor3}) represents the characteristic curves corresponding to invariants 
(\ref{firstcase}), while family (\ref{kor4}) corresponds to the characteristic invariants (\ref{secondcase}).}

\bigskip
\noindent
Let $(1,\varphi_0),\,(\varphi_0\neq 2\pi k,\,k \in \mathbb Z)$ be a point of the unit circumference. Then the following curves of families (\ref{curve1}) and (\ref{curve2}) pass respectively through this point:
\begin{equation} \label{kor6} \frac{2b(1-r\cos\vartheta)}{r^2-2r\cos\vartheta+1}-a^3\left\{ a^2 + \left[\frac{2br\sin\vartheta}{r^2-2r\cos\vartheta +1}+\frac{b\sin\varphi_0}{1-\cos\varphi_0}+a\sqrt{(a-b)/b}\right]^2\right\}^{-1} = 0,\;\;  \end{equation}
\begin{equation} \quad  0\leq\vartheta\leq\vartheta^*_1, \nonumber \end{equation}
\begin{equation} \label{kor7} \frac{2b(1-r\cos\vartheta)}{r^2-2r\cos\vartheta+1}-a^3\left\{ a^2 +
\left[\frac{2br\sin\vartheta}{r^2-2r\cos\vartheta +1}-\frac{b\sin\varphi_0}{1-\cos\varphi_0}-a\sqrt{(a-b)/b}\right]^2\right\}^{-1} = 0,  \end{equation}
\begin{equation} \vartheta^*_2 \leq\vartheta\leq2\pi, \nonumber \end{equation}
where $\vartheta^*_1, \vartheta^*_2$ are the polar angles of points of tangency of these curves with the circumference
(\ref{kor2}), respectively.

The relation
\begin{equation} \label{kor8} \frac{u_x(r,\vartheta,\varphi_0)}{u_y(r,\vartheta,\varphi_0)+1}=
\frac{f_1(\vartheta,\varphi_0)}{g_1(\vartheta,\varphi_0)},\end{equation}
is fulfilled for the characteristic curve (\ref{kor6}), and the relation
\begin{equation} \label{kor9} \frac{u_x(r,\vartheta,\varphi_0)}{u_y(r,\vartheta,\varphi_0)-1} =
\frac{f_2(\vartheta,\varphi_0)}{g_2(\vartheta,\varphi_0)},\end{equation}
is fulfilled along curve (\ref{kor7}), where
$$f_i(\vartheta,\varphi_0) = \frac{8a^3 bh_i r\sin\vartheta(1-\cos\vartheta)}{\left(a^2 + h_i^2\right)^2}+2b\;(r^2\cos2\vartheta-2r\cos\vartheta+1),$$
$$g_i(\vartheta,\varphi_0) = \frac{4a^3bh_i(r^2\cos2\vartheta - 2r\cos\vartheta+1)}{\left( a^2 + h_i^2\right)^2} - 4br\sin\vartheta (1-\cos\vartheta),   $$
$$h_i=\frac{2br\sin\vartheta}{r^2-2r\cos\vartheta+1}-\frac{b\sin\varphi_0}{1 -\cos\varphi_0}+(-1)^{i+1} a\sqrt{\frac{a-b}{b}}, \;i=1,2.$$
From (\ref{kor8}) and (\ref{kor9}) we conclude that the following equalities hold at the point $(1,\varphi_0)$:
\begin{equation} \label{kor10} \frac{u_x(1,\varphi_0)}{u_y(1,\varphi_0)+1}= \frac{f^0_1}{g^0_1},\end{equation}
\begin{equation} \label{kor11} \frac{u_x(1,\varphi_0)}{u_y(1,\varphi_0)-1}=\frac{f^0_2}{g^0_2},\end{equation}
where
$$f^0_i  =  (-1)^i \alpha \sin\varphi_0 + \cos\varphi_0,\,$$
$$g^0_i  =  (-1)^{i+1}  \alpha \cos\varphi_0 + \sin\varphi_0,$$
$$\alpha  =  \frac{2b^2}{a^2} \sqrt{\frac{a-b}{b}}, \;\,i=1,2.$$
Using relations (\ref{kor10}) and (\ref{kor11}) we define the values of the derivatives $u_x$ and $u_y$ at the point
$(1,\varphi_0)$
\begin{equation} u_x(1,\varphi_0)=\frac{2f^0_1f^0_2}{f^0_2g^0_1 - f^0_1g^0_2}, \nonumber \end{equation}
\begin{equation} u_y(1,\varphi_0)=\frac{f^0_1g^0_2 +f^0_2g^0_1}{f^0_2g^0_1 - f^0_1g^0_2}. \nonumber \end{equation}
Assuming that $\varphi_0$ takes all values on the circumference $x^2 + y^2 = 1,$ $0\leq\varphi_0\leq2\pi$, we have
\begin{equation} u_x(1,\vartheta) =\frac{\cos^2\vartheta - \frac{4b^3(a-b)}{a^4} \sin^2\vartheta}{\alpha}, \nonumber
\end{equation}
\begin{equation} u_y(1,\vartheta) = \frac{\alpha^2+1}{\alpha}\cos\vartheta\sin\vartheta. \nonumber
\end{equation}
Since the equalities
\begin{equation} u_\vartheta(1,\vartheta) = - u_x\sin\vartheta + u_y\cos\vartheta, \nonumber
\end{equation}
\begin{equation} u_r(1,\vartheta) = u_x\cos\vartheta + u_y\sin\vartheta, \nonumber
\end{equation}
are fulfilled, we obtain
\begin{equation} \label{kor11uj} u_\vartheta(1,\vartheta) = \alpha \sin\vartheta,
\end{equation}
\begin{equation} \label{kor12} u_r(1,\vartheta) = \frac{1}{\alpha}\cos\vartheta.
\end{equation}
Finally, integrating the value of (\ref{kor11uj}) on the circumference $x^2 + y^2 = 1, \;\; 0\leq\varphi_0\leq2\pi$, we define also the
function
\begin{equation} \label{kor13} u(1,\vartheta) = -\alpha\cos\vartheta+c. \end{equation}

If we consider the initial problem: find a function $u(x,y)$,  which satisfies equation (\ref{equ2}) with initial data defined by functions (\ref{kor12}),
(\ref{kor13}) on the unit circumference $\gamma $: $r=1$, $x=r \cos \vartheta$, $y=r \sin \vartheta$, then the definition domain of a solution of this problem is constructed by the set of characteristic curves represented by formulas (\ref{kor3}), (\ref{kor4}). By the structure of these characteristic curves we conclude that for $a>b$ the definition domain of a solution of the initial problem is the domain of the plane $(x,y)$ which lies to the left of the straight line
$x=1$, with a gap given the circumference (\ref{kor2}).

\end{document}